\documentclass{lms}

\usepackage{amsmath,amssymb,amsfonts}

\usepackage{xspace,enumerate}
\usepackage{booktabs}
\usepackage{color,epsfig,graphics}
\newtheorem{scase}{Case}[section]

\usepackage{calc}
\usepackage{burdges-JLMS}

\newcommand\Uir{U_{0,r}}

\newcommand\reB{\tilde{B}}
\newcommand\reU{\tilde{U}}
\newcommand\reE{\tilde{E}}
\newcommand\reF{\tilde{F}}

\newcommand\components{\mathcal{E}}
\newcommand\recomponents{\tilde{\components}}

\title[A new trichotomy theorem]%
{A new trichotomy theorem for \\ groups of finite Morley rank}
\author{Alexandre Borovik \and Jeffrey Burdges}
\extraline{Burdges was supported by NSF grant DMS-0100794, and
 Deutsche Forschungsgemeinschaft grant Te 242/3-1.}
\classno{03C60 (primary), 20G99 (secondary)}

\begin{document}
\maketitle

\begin{abstract}
There is a longstanding conjecture, due to Gregory Cherlin and Boris Zilber,
that all simple groups of finite Morley rank are simple algebraic groups.
Here we will conclude that a simple $K^*$-group of finite {M}orley rank
 and odd type either has normal rank at most two, or else
 is an algebraic group over algebraically closed field of characteristic not 2.
To this end, it suffices, by \cite{Bu07a,BBN04,BCJ}, to produce a 
proper 2-generated core in groups with \Prufer rank two and
normal rank at least three; which is what is proven here.
Our final conclusion constrains the Sylow 2-subgroups available to
a minimal counterexample, and finally proves the trichotomy
theorem in the nontame context.
\end{abstract}

The Algebraicity Conjecture for simple groups of finite Morley rank,
also known as the Cherlin-Zilber conjecture, states that simple groups
of finite Morley rank are simple algebraic groups over algebraically
closed fields.  In the last 15 years, considerable progress has been
made by transferring methods from finite group theory; however,
the conjecture itself remains decidedly open.  In the formulation of
this approach, groups of finite Morley rank are divided into four
types, odd, even, mixed, and degenerate, according to the structure
of their Sylow 2-subgroup.  For even and mixed type the Algebraicity
Conjecture has been proven, and connected degenerate type groups
are now known to have trivial Sylow 2-subgroups \cite{BBC}.
Here we concern ourselves with the ongoing program to analyze a
minimal counterexample to the conjecture in {\em odd type}, where
the Sylow 2-subgroup is divisible-abelian-by-finite.

The present paper lies between the high \Prufer rank, or generic,
where general methods are used heavily, and the ``end game''
where general methods give way to consideration of special cases.
In the first part, the Generic Trichotomy Theorem \cite{Bu07a} says
that a {\em minimal} non-algebraic simple group of finite Morley rank
has \Prufer rank at most two.  Thus we may consider small groups
whose simple sections are restricted to $\PSp_4$, $\G_2$, $\PSL_3$,
and $\PSL_2$.  In the next stage, we hope to proceed via the analysis
of components in the centralizers of {\em toral} involutions; however,
the existence of such a component requires further constraints on the
Sylow 2-subgroup \cite[p.~196]{BuPhd}.
In particular, the normal 2-rank provides another analog of Lie rank,
more delicate and traditional than the \Prufer 2-rank,
 which needs to be controlled.

Therefore the present paper argues that
\begin{namedtheorem}{Theorem}
Any simple $K^*$-group $G$ of finite Morley rank and odd type,
 with \Prufer 2-rank at most two and normal 2-rank at least 3,
has a proper 2-generated core, i.e.\ $\Gamma^0_{S,2}(G) < G$.
\end{namedtheorem}

For tradition and technical reasons, we prefer to record this result as
a tameness free version of the old trichotomy theorem from \cite{Bo95},
using the Generic Trichotomy Theorem \cite{Bu07a} (page \pageref{GIT})
when the \Prufer 2-rank is at least 3.

\begin{namedtheorem}{Trichotomy Theorem}
Let $G$ be a simple $K^*$-group of finite Morley rank and odd type.
Suppose either that $G$ has normal 2-rank $\geq3$, or else
 has an eight-group centralizing the 2-torus $S^\o$.
Then one of the following holds.
\begin{conclusions}
\item $G$ has a proper weak 2-generated core, i.e.\ $\Gamma^0_{S,2}(G) < G$.
\item $G$ is an algebraic group over an algebraically closed field of characteristic \noteql2.
\end{conclusions}
\end{namedtheorem}

It follows immediately from normal 2-rank $\geq 3$ that some
eight-group $\cong (\Z/2\Z)^3$ centralizes some maximal 2-torus.
So we work throughout with this more general condition, which
we plan to exploit in the later stages described above.
By Fact \ref{Asch46.2-3} below, normal 2-rank $\geq 3$ also implies
that the 2-generated core $\Gamma_{S,2}(G)$ from \cite{Bo95} and
our weak 2-generated core $\Gamma^0_{S,2}(G)$ coincide.


We summarize the present status of the classification as follows.

\begin{namedtheorem}{Status}\label{status}
Let $G$ be a simple $K^*$-group of finite Morley rank and odd type,
 with normal 2-rank $\geq3$ and \Prufer 2-rank $\geq2$.
Then $G$ is an algebraic group over an algebraically closed field of characteristic \noteql2.
\end{namedtheorem}

This result follows immediately from two previous papers,
 our Trichotomy Theorem, and Fact \ref{Asch46.2-3} below.

\begin{namedtheorem}{Strong Embedding Theorem}[{\cite{BBN04}}]
Let $G$ be a simple $K^*$-group of finite Morley rank and odd type,
 with normal 2-rank $\geq3$ and \Prufer 2-rank $\geq2$.
Suppose that $G$ has a proper 2-generated core $M = \Gamma_{S,2}(G) < G$.
Then $G$ is a minimal connected simple group, and $M$ is strongly embedded.
\end{namedtheorem}

\begin{namedtheorem}{Minimal Simple Theorem}[{\cite{BCJ}}]
Let $G$ be a minimal connected simple group of finite Morley rank and of odd type.
Suppose that $G$ contains a proper definable strongly embedded subgroup $M$.
Then $G$ has \Prufer 2-rank one.
\end{namedtheorem}

The argument presented below proceeds via the analysis of components of
the centralizers of involutions.  A purely model theoretic version has recently
been obtained by Burdges and Cherlin.

\section{Cores and components}

We begin by recalling the various consequences of the absences of
 a proper 2-generated core, as laid out in \cite[\S2]{Bu07a}.
Of primary importance is the role of quasisimple components
 of the centralizers of toral involutions.


Much of our time will be spent analyzing a so-called simple $K^*$-group
of finite Morley rank.  A $K^*$-group is a group whose {\em proper}
definable simple sections are all algebraic.  $K^*$-groups are analyzed by
examining various proper subgroups, especially the centralizers of involutions.
A group is said to be a $K$-group if {\em all} definable simple sections are
algebraic, and this property is holds for proper subgroups of our $K^*$-group.

An algebraic group is said to be {\em reductive} if it has no unipotent
radical, and a reductive group is a central product of semisimple
algebraic groups and algebraic tori.  In a simple (even reductive)
algebraic group, over a field of  characteristic \noteql2, the centralizer
of an involution is itself reductive.  In this section, we establish,
in the absence of a proper 2-generated core, that the centralizers of
involutions in our simple $K^*$-group are ``somewhat reductive''.

\subsection{Proper 2-generated cores}

\begin{definition}
Consider a group $G$ of finite Morley rank and a 2-subgroup $S$ of $G$.
We define the {\em 2-generated core} $\Gamma_{S,2}(G)$ of $G$
(associated to $S$) to be the definable hull of the group generated by
all normalizers of four-subgroups in $S$.
We also define the {\em weak 2-generated core} $\Gamma^0_{S,2}(G)$
of $G$ (associated to $S$) to be the definable hull of all normalizers of
four-subgroups $U \leq S$ with $m(C_S(U)) > 2$.
We say that $G$ has a {\em proper 2-generated core} or
 a {\em proper weak 2-generated core} when, for a Sylow 2-subgroup $S$,
 $\Gamma_{S,2}(G) < G$ or $\Gamma^0_{S,2}(G) < G$, respectively.
\end{definition}

Both notions are well-defined because the Sylow 2-subgroups of
a group of finite Morley rank are conjugate \cite[\qTheorem 10.11]{BP,BN}.

\begin{fact}[{\cite[\qFact 1.20-2]{Bu07a}; compare \cite[46.2]{Asch}}]%
\label{Asch46.2-3}
Let $G$ be a group of finite Morley rank, and let $S$ be a 2-subgroup of $G$.
If $S$ has normal 2-rank $\geq 3$, then $\Gamma^0_{S,2}(G) = \Gamma_{S,2}(G)$.
\end{fact}

For an elementary abelian 2-group $V$ acting definably on $G$,
we define $\Gamma_V(G)$ to be the group generated by the
connected components of centralizers of involutions in $V$.
$$ \Gamma_V(G) = \gen{C^\o_G(v) : v\in V^\#}\mathperiod $$

Out most basic tool for producing a proper 2-generated core is the following.

\begin{fact}[{\cite[\qProposition 1.22]{Bu07a}}]\label{GITcondA}
Let $G$ be a simple $K^*$-group of finite Morley rank and odd type,
with $m(G)\geq3$, and let $S$ be a Sylow 2-subgroup of $G$.
Suppose that $\Gamma_E(G) < G$ for some
 four-group $E \leq G$ with $m(C_G(E)) > 2$.
Then $G$ has a proper weak 2-generated core. 
\end{fact}

\subsection{Components and descent}

\begin{definition}
A quasisimple subnormal subgroup of a group $G$ is called a
{\em component} of $G$ (see \cite[p.~118 ii]{BN}).
We define $E(G)$ to be the connected part of the product of components
of $G$, or equivalently the product of the components of $G^\o$
(see \cite[\qLemma 7.10iv]{BN}).
Such components are normal in $G^\o$ by \cite[\qLemma 7.1iii]{BN},
 and indeed $E(G) \normal G$.
\end{definition}

A connected reductive algebraic group $G$ is a central product
of $E(G)$ and an algebraic torus.

We require several plausible ``unipotent radicals''
 to define our notion of partial reductivity.

\begin{definition}
The {\em Fitting subgroup} $F(G)$ of a group $G$
 is the subgroup generated by all its nilpotent normal subgroups.
\end{definition}

In any group of finite Morley rank, the Fitting subgroup is
 itself nilpotent and definable \cite[Theorem 7.3]{Bel87,Ne91,BN},
and serves as a notion of unipotence in some contexts.

\begin{definition}
A connected definable $p$-subgroup of bounded exponent inside
a group $H$ of finite Morley rank is said to be {\it $p$-unipotent}.
We write $U_p(H)$ for the subgroup generated by
 all $p$-unipotent subgroups of $H$.
\end{definition}

If $H$ is solvable, then $U_p(H) \leq F^\o(H)$ is $p$-unipotent itself
(see \cite[\qCorollary 2.16]{CJ01} and \cite[\qFact 2.36]{ABC97}).

\begin{definition}
We say that a connected abelian group of finite Morley rank is
{\it indecomposable} if it has a unique maximal proper
definable connected subgroup, denoted $J(A)$ (see \cite[\qLemma 2.4]{Bu03}).
We define the {\it reduced rank} $\rr(A)$ of
a definable indecomposable abelian group $A$
to be the Morley rank of the quotient $A/J(A)$,
i.e.\ $\rr(A) = \rk(A/J(A))$.
For a group $G$ of finite Morley rank, and any integer $r$, we define
$$ \Uir(G) = \Genst{A \leq G}{%
\parbox{\widthof{$A$ is a definable indecomposable group,}}%
{$A$ is a definable indecomposable group, \\
\hspace*{10pt} $\rr(A) = r$, and $A/J(A)$ is torsion-free}}\mathperiod $$
We say that $G$ is a {\it $\Uir$-group} if $U_{0,r}(G)=G$, and
 set $\rr_0(G) = \max \{r \mid \Uir(G) \neq 1 \}$.
\end{definition}

We view the reduced rank parameter $r$ as a {\em scale of unipotence},
with larger values being more unipotent.  By \cite[\qTheorem 2.16]{Bu03},
the ``most unipotent'' groups, in this scale, are nilpotent.

\begin{definition}
In a group $H$ of finite Morley rank,
 we write $O(H)$ for the subgroup generate by
 the definable connected normal subgroups without involutions.
\end{definition}

If $H$ is a $K$-group of finite Morley rank, then $O(H)$ is solvable,
as simple algebraic groups always contain 2-torsion.

Our approach to reductivity begins with the following fact.

\begin{fact}[{\cite[\qTheorem 5.12]{Bo95}}]\label{KgrpO}
Let $G$ be a connected $K$-group of finite Morley rank and odd type
with $O(G) = 1$.  Then $G = F^\o(G) * E(G)$ is isomorphic to
a central product of quasisimple algebraic groups over
algebraically closed fields of characteristic \noteql2 and
of a definable normal divisible abelian group $F^\o(G)$.
\end{fact}

\noindent However, a more subtle definition is required to find an
applicable version of this fact.  The following definition was applied
in \cite{Bu07a}, under the assumption of \Prufer rank $\geq 3$.

\begin{definition}\label{def:rbarstar}
Consider a simple group $G$ of finite Morley rank and
let $X$ be a subgroup of $G$ with $m(X)\geq3$.
We write $I^0(X) := \{ i\in I(G) : m(C_X(i))\geq3 \}$
for the set of involutions from eight-groups in $X$.
We define $\rr^*(X)$ to be the supremum of $\rr_0(k^*)$ as $k$ ranges
over the base fields of the algebraic components of the quotients
$C^\o_G(i) / O(C_G(i))$ associated to involutions $i\in I^0(X)$.
\end{definition}

Clearly $\rr^*(G)$ is the maximum of $\rr^*(E)$
 as $E$ ranges over eight-groups in $G$.
We recall that, for a nonsolvable group $L$ of finite Morley rank,
$\Uir(L)$ and $U_p(L)$ need not be solvable, as quasisimple algebraic
groups are generated by the unipotent radicals of their Borel subgroups.
We exploit this in the following central definition.

\begin{definition}\label{def:reE}
We continue in the notation of Definition \ref{def:rbarstar}.
For a definable subgroup $H$ of $G$, we define $\reU_X(H)$ to be
 the subgroup of $H$ generated by $U_p(H)$ for $p$ prime
 as well as by $\Uir(H)$ for $r > \rr^*(X)$.
As an abbreviation, we use $\reF_X(H)$ to denote $F^\o(\reU_X(H))$,
and $\reE_X(H)$ to denote $E(\reU_X(H))$.
We use $\recomponents^X_Y$ to denote the set of components
of $\reE_X(C_G(i)) = E(\reU_X(C_G(i)))$ for $i\in I^0(Y)$ with $Y \leq X$,
 and set $\recomponents_X = \recomponents^X_X$.
\end{definition}

$\reU_X(H)$ is the subgroup of $H$ which is generated by its
unmistakably unipotent subgroups.  These definitions are all
sensitive to the choice of $X$, which is usually a fixed eight-group.

\begin{definition}\label{def:reB}
We say that a simple $K^*$-group $G$ with $m(G)\geq3$ satisfies
 the {\em $\reB$-property} if, for every 2-subgroup $X \leq G$ with $m(X)\geq3$
 and every $t\in I^0(X)$, the group $\reU_X(O(C_G(t)))$ is trivial.
\end{definition}

We recall from \cite{Bu07a} that the $\reB$-property holds
 in the absence of a proper 2-generated core.
The proof of this result explains the definition of $\rr^*$.

\begin{fact}[{\cite[\qTheorem 2.9]{Bu07a}}]\label{rebalancing}
Let $G$ be a simple $K^*$-group of finite Morley rank and odd type
with $m(G) \geq 3$.  Then either
\begin{conclusions}
\item $G$ has a proper weak 2-generated core,  or
\item $G$ satisfies the $\reB$-property.
\end{conclusions}
\end{fact}

For us, the point of the $\reB$-property is that it ensures the existence of a
well behaved family of components in the centralizers of involutions.

\begin{definition}
Given a simple group $G$ of finite Morley rank and a set of involutions $J$,
we say a family of components $\components$ from the centralizers
of involutions in $J$ is {\em descent inducing for $J$} if
\begin{quote}
For every $K\in \components$
 and every involution $i\in J$ which normalizes $K$,
there are components $L_1,\ldots,L_n \in \components$
 with $L_k \normal C^\o_G(i)$ such that $E(C_K(i)) \leq L_1 * \cdots * L_k$.
\end{quote}
In this situation we say that $K$, or just $E(C_K(i))$, descends.
\end{definition}

The $\reB$-property provides us with such a family of components.

\begin{lemma}\label{KgrpO_reU_descent}
Let $G$ be a simple $K^*$-group of finite Morley rank and odd type
with $m(G)\geq3$ which satisfies the $\reB$-property.
Then, for every 2-subgroup $X \leq G$ with $m(X)\geq3$, the family
  of components $\recomponents_X$ is descent inducing for $I^0(X)$.
\end{lemma}

We will extract this result from the following facts.

\begin{fact}[{\cite[\qCorollary 2.11]{Bu07a}}]\label{KgrpO_reU}
Let $G$ be a simple $K^*$-group of finite Morley rank and odd type
with $m(G)\geq3$ which satisfies the $\reB$-property.
Then, for every 2-subgroup $X \leq G$ with $m(X)\geq3$ and every $i\in I^0(X)$,
we have $\reU_X(C_G(i)) = \reE_X(C_G(i)) * \reF_X(C_G(i))$
 and $\reF_X(C_G(i))$ is abelian.
\end{fact}

\begin{fact}[{\cite[\qProposition 2.13]{Bu07a}}]\label{control_field}
Let $H$ be a groups of finite Morley rank which is isomorphic to
 a linear algebraic groups over an algebraically closed field $k$.
Then
\begin{conclusions}
\item If $\Uir(H) \neq 1$ for some $r > \rr_0(k^*)$
   then $\Char(k) = 0$ and $\rk(k) = r$.
\item If $U_p(H) \neq 1$ then $\Char(k) = p$.
\end{conclusions}
If $H$ is quasisimple, these conditions imply
 $U_p(H) = H$ and $\Uir(H) = H$, respectively.
\end{fact}


\begin{fact}[{\cite[\qTheorem 8.1]{St2}; \cite[\qFact 1.6]{Bu07a}}]\label{Creductive2}
Let $G$ be a quasisimple algebraic group over an algebraically closed field.
Let $\phi$ be an algebraic automorphism of $G$ whose order is finite and
relatively prime to the characteristic of the field.
Then $C^\o_G(\phi)$ is nontrivial and reductive.
\end{fact}



\begin{proof}[of Lemma \ref{KgrpO_reU_descent}]
Consider a component $K\in \recomponents_X$
 and an involution $i\in I^0(X)$ which normalizes $K$.
By Fact \ref{Creductive2}, $C_K^\circ(i)$ is reductive.
By Fact \ref{control_field}
 $\reU_X(C_G(i)) \geq \reE_X(C_K(i)) = E(C_K(i))$.
As $E(C_K(i))$ is nonabelian, Fact \ref{KgrpO_reU} yields
 a set of components $L_1,\ldots,L_n \in \recomponents_X$
 with $L_k \normal C^\o_G(i)$ such that
\[ E(C_K(i)) \leq L_1 * \cdots * L_n \mathperiod \]
\end{proof}

Our descent inducing family of components generates $G$.

\begin{fact}[{\cite[\qTheorem 2.18]{Bu07a}}]\label{rebalanced_Asar2}
Let $G$ be a $K^*$-group of finite Morley rank and odd type with $m(G)\geq3$.
Suppose that there is a four-group $E \leq G$
 which centralizes a Sylow\o 2-subgroup $T$ of $G$, and that
 there is an eight-group $X$ in $C_G(T)$ which contains $E$.
Then either
\begin{conclusions}
\item $G$ has a proper weak 2-generated core, or else
\item $\gen{\recomponents^X_E} = \gen{\reE_X(C_G(z)) : z\in E^\#} = G$.
\end{conclusions}
\end{fact}

\subsection{Automorphisms}

The following two facts ensure that involutions acting upon
quasisimple components are understood.

\begin{definition}
Given an algebraic group $G$, a maximal torus $T$ of $G$, and a
Borel subgroup $B$ of $G$ which contains $T$, we define the group
$\Gamma$ of {\em graph automorphisms} associated to $T$ and $B$,
to be the group of algebraic automorphisms of $G$ which normalize
both $T$ and $B$.
\end{definition}


\begin{fact}[{\cite[\qTheorem 8.4]{BN}}] \label{autalg}
Let $G \rtimes H$ be a group of finite Morley rank,
 with $G$ and $H$ definable.
Suppose that $G$ is a quasisimple algebraic group over an
 algebraically closed field, and $C_H(G)$ is trivial.
Let $T$ be a maximal torus of $G$ and
 let $B$ be a Borel subgroup of $G$ which contains $T$.
Then, viewing $H$ as a subgroup of $\Aut(G)$, we have $H\leq \Inn(G)\Gamma$,
where $\Inn(G)$ is the group of inner automorphisms of $G$ and $\Gamma$
is the group of graph automorphisms of $G$ associated to $T$ and $B$.
\end{fact}

% ...

\begin{fact}
\label{invautalg}
Let $G$ be an infinite quasisimple algebraic group over an algebraically
closed field of characteristic \noteql2 and let $\phi$ be an involutive
automorphism centralizing some Sylow\o 2-subgroup $T$ of $G$.
Then $\phi$ is an inner automorphism induced by some element of $T$.
\end{fact}

%

\begin{proof} 
We observe that $H := C_G(T)$ is a maximal torus of $G$ containing $T$.
Let $B$ be a Borel subgroup containing $H$.
By Fact \ref{autalg}, $\phi = \alpha \circ \gamma$ where
 $\gamma$ is a graph automorphism normalizing $H$ and $B$, and
 $\alpha(x) = x^t$ is an inner automorphism induced by some element $t\in G$.
Since $\phi$ and $\gamma$ normalize the $H$,
 we know that $\alpha$ does as well, so $t\in N_G(H)$.

Following \cite[p.~17--18]{Carter93}, we consider the root system
$\Phi$ to be a subset of the set of algebraic homomorphisms $\Hom(H,k^*)$
from $H$ to the multiplicative group of the field $k^*$.
The action of both the Weyl group and the graph automorphisms of $G$ on
the root system is the natural action on $\Hom(H,k^*)$.
If $t \notin H$, then $t$ is a representative for a nontrivial Weyl group
element, would not preserve the set of positive roots in $\Phi$.
On the other hand, $\phi$ centralizes $T$,
 so it preserves the set of positive root of $\Phi$.
Since $\gamma$ normalizes $B$, it preserves the set of positive roots too.
So $t\in H$.
Now $\con\phi=\gamma(H)$ is a graph automorphism itself.
By \cite[Table 4.3.1\ p.~145]{GLS3},
 a nontrivial involutive graph automorphism of
 a quasisimple algebraic group $G$ never centralizes a maximal torus of $G$.
So $\phi$ must be an inner automorphism.
\end{proof}

\smallskip

Table \ref{T:alggrpdata} contains necessary information about 
 conjugacy classes of involutions and their centralizers,
 in Lie rank two quasi-simple groups
(see \cite[Table 4.3.1\ p.\ 145 \& Table 4.3.3\ p.\ 151]{GLS3}).

\begin{table}[h]
\begin{center}
\begin{tabular}{llllll} \toprule
$G$ & $\bar\Gamma$ & $Z$ & $i$ & $C^\o_G(i)$ \\
\midrule
$\SL_2$ & 1 & $\Z/2\Z$ & inner & $k^*$ \\
$\PSL_2$  & 1& 1 & inner & $k^*$ \\
$\SL_3$ & $\Z/2\Z$ & $\Z/3\Z$
   & inner & $\SL_2 * k^*$ \\
& & & graph & $\PSL_2$ \\
$\PSL_3$ & $\Z/2\Z$ & 1 
   & inner & $\SL_2 * k^*$ \\
& & & graph & $\PSL_2$ \\
$\Sp_4$ & 1 & $\Z/2\Z$
   & inner & $\SL_2 \times \SL_2$ \\
& & & inner & $\SL_2 * k^*$ \\
$\PSp_4$ & 1 & 1
   & inner & $\SL_2 * \SL_2$ \\
& & & inner & $\PSL_2 \times k^*$ \\
$\G_2$ & 1 & 1 & inner & $\SL_2 * \SL_2$ \\ \bottomrule
\end{tabular}
\end{center}
\caption{Data on Chevalley Groups}\label{T:alggrpdata}%
\label{invSL3}\label{invSp4}\label{invG2}\label{findSp4}
\end{table}

\section{Proof of the Trichotomy Theorem}\label{sec:Trico}

We prove our main result in this section.

\begin{namedtheorem}{Trichotomy Theorem}
Let $G$ be a simple $K^*$-group of finite Morley rank and odd type.
Suppose either that $G$ has normal 2-rank $\geq3$, or else
 has an eight-group centralizing the 2-torus $S^\o$.
Then one of the following holds.
\begin{conclusions}
\item $G$ has a proper 2-generated core.
\item $G$ is an algebraic group over an algebraically closed field of characteristic \noteql2.
\end{conclusions}
\end{namedtheorem}

In high \Prufer rank, \cite{Bu07a} produces an algebraic group.
In other situations, we will reach a contradiction below.

\subsection{Consequences of $\Gamma^0_{S,2}(G) = G$.}

We observe that an elementary abelian 2-group which is normal in
a Sylow 2-subgroup $S$ is centralized by $S^\o$.  So we may begin
proof of the trichotomy theorem by assuming the following.

\begin{hypothesis}
Let $G$ be a simple $K^*$-group of finite Morley rank and odd type
which satisfies the following.
\begin{hypotheses}
\item  $G$ has an eight-group centralizing the 2-torus $S^\o$.
\item  $G$ does not have a proper weak 2-generated core, i.e. $\Gamma^0_{S,2}(G) = G$.
\end{hypotheses}
\end{hypothesis}

Let $S$ be a Sylow 2-subgroup of $G$, and let $A \leq S$ be an
eight-group such that $[A,S^\o] = 1$.
We may also assume $\Omega_1(S^\o) \leq A$.
The high \Prufer rank case is covered by the following.

\begin{namedtheorem}{Generic Trichotomy Theorem}[{\cite{Bu07a}}]\label{GIT}
Let $G$ be a simple $K^*$-group of finite Morley rank and odd type
with \Prufer 2-rank $\geq 3$.  Then either
\begin{conclusions}
\item $G$ has a proper weak 2-generated core, or
\item $G$ is an algebraic group over an algebraically closed field of characteristic \noteql2.
\end{conclusions}
\end{namedtheorem}

As we have assumed that $\Gamma^0_{S,2}(G) = G$,
 the Generic Trichotomy Theorem allows us to assume that
\begin{center} $G$ has \Prufer 2-rank 1 or 2. \end{center}
\smallskip
By Fact \ref{rebalancing}, $G$ satisfies the $\reB$-property too.
In particular, the family of components $\recomponents_A$
 is descent inducing for $A^\#$ by Lemma \ref{KgrpO_reU_descent}.
We will make extensive use of this family.

\begin{definition}
Let $\components := \recomponents_A$, and
let $\components^* \subset \components$ be the set of those
components in $\components$ with Lie rank two.
\end{definition}

Since $m(A) \geq 2$ and $[A,S^\o] = 1$,
 we have $\gen{\recomponents^A_{\Omega_1(S^\o)}} = G$
 by Fact \ref{rebalanced_Asar2}.
In particular, $\recomponents^A_{\Omega_1(S^\o)} \neq \emptyset$.
As $G$ is simple, we also obtain the following.
\begin{texteqn}{\tag{$\star$}}
For $M\in \components$,
  there is a $v\in \Omega_1(S^\o)^\#$ such that $E(C_G(v)) \not\leq M$.
\end{texteqn}

We next show that $G$ has a Lie rank two component,
 i.e. $\components^* \neq \emptyset$.
This contains our final application of $\Gamma^0_{S,2}(G) = G$.
It also shows that $\pr(G) \neq 1$, and hence that $\pr(G)=2$.

\begin{lemma}\label{Lieranktwocomp_exists}
For any component $K\in \components$,
there is some component $M\in \components^*$ containing $K$.
\end{lemma}

We will need the following fact.

\begin{fact}[{\cite[\qLemma 8.1]{Bo95}; compare \cite[Proposition I.1.1]{WalterClsInv}}]%
\label{f:newtrico_walter}
Let $H$ be a connected $K$-group of finite Morley rank with $H = F(H) E(H)$.
Let $t$ be a definable involutive automorphism of $H$ and
let $L$ be a component of $E(C_H(t))$.  Then there is a component
$K \normal E(H)$ such that one of the following holds.
\begin{conclusions}
\item $K = K^t$ and $L \normal E(C_K(t))$.
\item $K \neq K^t$ and $L \normal E(C_{K^t K}(t))$.
\end{conclusions}
\end{fact}

\begin{proof}
We may assume $H = E(H)$ as $L \leq E(H)$.
For any component $K$ of $E(H)$, either $K^t = K$, or else
 $C_{K^t K}(t)$ is contained in the diagonal $\Delta(K)$ of
 the central product $K^t * K$.
We take $\Delta(K) := K$ in the first situation, to simplify notation.
Now $C_H(t) \leq \Delta(K_1) * \cdots * \Delta(K_n)$.
As all $\Delta(K_i)$s are fixed by $i$,
 $C_H(t) = C_{\Delta(K_1)}(t) * \cdots * C_{\Delta(K_n)}(t)$ too.
So there is a unique $\Delta(K_i)$ containing $L$, as desired.
\end{proof}

\smallskip


\begin{proof}[Proof of Lemma \ref{Lieranktwocomp_exists}]
We may assume that $K\notin \components^*$ has type $\pPSL_2$.
Since $A$ centralizes $S^\o$, $S^\o$ is a Sylow\o 2-subgroup of $C^\o_G(y)$
 for any $y\in A^\#$, and $A$ fixes all components of $E(C_G(y))$.
In particular, $A$ normalizes $K$.

Since $K$ has no graph automorphisms, 
 $A$ acts on $K$ via inner automorphisms by Fact \ref{autalg},
 and these inner automorphisms normalize the 2-torus $K \cap S^\o$.
Since a four-group in $\Aut(K) \cong \PSL_2$ normalizing this 2-torus 
 has an involution inverting the torus and $A$ centralizes the 2-torus,
 we have that $m(C_A(K)) \geq m(A)-1 \geq 2$.
Since $K\in \components$ induces descent for $A^\#$,
 $K \leq E(C_G(x))$ for all $x\in C_A(K)$.

Now suppose that $K \normal E(C_G(x))$ for all $x\in C_A(K)$.
So $K \normal \Gamma_{C_A(K)}(G)$ and hence $\Gamma_{C_A(K)}(G) < G$.
By Fact \ref{GITcondA},
 $G$ must have a proper weak 2-generated core.
Thus
$$ \textrm{There is an $x\in C_A(K)$ such that
  $K$ is not normal in $E(C_G(x))$.} $$

Fix such an $x$.
Suppose toward a contradiction that $E(C_G(x))$ is not quasisimple.
Since $\pr(G) \leq 2$, we have
 $E(C_G(x)) = L_1 * L_2$ with $L_i$ of $\pPSL_2$ type.
Since $K \normal E(C_G(y))$ for some $y\in A^\#$,
 we have $K \normal E(C_{L_1*L_2}(y))$ too.
By Fact \ref{f:newtrico_walter},
 there is a component $L \normal E(L_1*L_2) = L_1 * L_2$ such that either
\begin{enumerate}
\item $L^y = L$ and $K \normal E(C_L(y))$, or
\item $L^y \neq L$ and $K \normal E(C_{L^y L}(y))$
\end{enumerate}
Since $L \cong \pPSL_2$,
 $y$ acts on $L$ by an inner automorphism by Fact \ref{autalg},
so $C_L(y)$ is an algebraic subgroup of $L$.
In case (i), either $C_L(y) = L$ or $C _L(y)$ is solvable,
 and both are contradictions.
In case (ii), we have $L^y \neq L$, contradicting $[A,S^\o] = 1$.
\end{proof}

Since $\components \neq \emptyset$,
 we now have $\components^* \neq \emptyset$,
 and thus $$ \pr(G)=2\mathperiod $$

\subsection{Lie rank two components}

A brief inspection of \cite[Table 4.3.1\ p.~145]{GLS3}
will reveal that quasisimple algebraic groups almost never have
centralizers of involutions which are themselves quasisimple of
the same Lie rank, the only exceptions being $B_4$ inside $F_4$,
$A_7$ inside $E_7$, and $D_8$ inside $E_8$.  In particular,
the Lie rank two quasisimple algebraic groups $\pPSL_3$,
$\pPSp_4$, and $G_2$ never have Lie rank two quasisimple
centralizers.
We view this observation as inspirational, and
 proceed to reach a contradiction by considering the various possible
 types of components in $\components^*$.
Our hypotheses are as follows,
 given the analysis of the preceding subsection.

\begin{hypothesis}
Let $G$ be a simple $K^*$-group of finite Morley rank and odd type
which satisfies the following.
\begin{hypotheses}
\item $G$ has an eight-group $A$ centralizing the 2-torus $S^\o$.
\item $G$ has a nonempty family of components $\components$,
 from the centralizers of involutions in $A^\#$,
 which is descent inducing for $A^\#$.
\item The set $\components^* \subset \components$ of
 components with Lie rank two is nonempty.
\item For all $M\in \components$,
  there is a $v\in \Omega_1(S^\o)^\#$ such that $E(C_G(v)) \not\leq M$.
\end{hypotheses}
\end{hypothesis}

Fix some $M \in \components^*$. 
As $M$ must be isomorphic to one of $\pPSL_3$, $\pPSp_4$, or $G_2$,
 we proceed by analyzing each of these cases separately.
Each case will either reach a contradiction, or arrive at a previous treated case.

\begin{scase}\label{Lieranktwocomp_Sp4}
$\components^*$ contains a component $M \cong \Sp_4$.
\end{scase}

\begin{analysis}
Since $S^\o$ is a Sylow\o 2-subgroup of $M$,
 and the central involution of $\Sp_4$ is toral,
there is an $x\in \Omega_1(S^\o)$ with $M \leq E(C_G(x))$,
 and $x$ is the central involution of $M$.
It follows that the group $E(C_G(x))$ is quasisimple,
 so $M = E(C_G(x)) \cong \Sp_4$.
Since $M \not\leq C_M(u)$ for some $u \in \Omega_1(S^\o) \setminus \gen{x}$,
 $u$ is not conjugate to $x$ in $M$.
Since $\Sp_4$ has exactly two conjugacy classes of involutions
by Table \ref{invSp4}, the two distinct involutions
$u,v \in \Omega_1(S^\o) \setminus \gen{x}$ are conjugate.
By Table \ref{invSp4}, $C^\o_M(u) \cong \SL_2 \times \SL_2$.

By $(\star)$ and the conjugacy of $u$ and $v$,
 we have $E(C_G(u)) \not\leq M$.
Since $u\in S^\o \leq C^\o_M(u)$
 and $C^\o_M(u) \cong \SL_2 \times \SL_2$,
 we find $u\in Z(E(C_M(u)))$.
By Table \ref{findSp4}, $K_u := E(C_G(u)) \cong \Sp_4$.
Since $u$ and $v$ are conjugate, $K_v := E(C_G(v)) \cong \Sp_4$ as well.
Since $M$ is descent inducing for $A^\#$,
 we have $K_u,K_v \in \components^*$.

Since $[A,S^\o]=1$, the group $A$ acts on $M$
 by inner automorphisms induced by elements of $S^\o$
 by Fact \ref{invautalg}.
So there is an involution $w\in A \setminus \Omega_1(S^\o)$
 such that $M \leq E(C_G(w))$.
By Table \ref{invSL3}, 
 only $\Sp_4$ has Lie rank two and
 $\SL_2 \times \SL_2$ inside the centralizer of an involution,
so $E(C_G(w)) \cong \Sp_4$.  Thus $x \in Z(E(C_G(w)))$.
Since $M$ induces descent for $A^\#$, we have $E(C_G(w)) \in \components$.
Since $E(C_G(w))$ induces descent for $A^\#$,
 we have $E(C_G(w)) \leq E(C_G(x)) = M$.
So $E(C_G(w)) = E(C_G(x)) = M$ and $E(C_G(w)) = M$.

Since $w$ centralizes $u$, the involution $w$ normalizes $K_u$.
Since $K_u$ is descent inducing for $A^\#$,
we have $E(C_{K_u}(w)) \leq M \cap K_u$
 and $E(C_{K_u}(w)) \geq C^\o_M(u) \cong \SL_2 \times \SL_2$
 by Table \ref{invSp4}.
Since $w$ centralizes $S^\o$, Fact \ref{invautalg} says that
 $w$ acts by an inner automorphism on $K_u$.
Since $K_u \not\leq M = E(C_G(w))$,
 we know that $w$ does not centralize $K_u$.
By Table \ref{invSp4},
 $C^\o_{K_u}(w) \cong \SL_2 \times \SL_2$ or $\SL_2 * k^*$,
 so $C^\o_{K_u}(w) \cong \SL_2 \times \SL_2$,
 and $C^\o_{K_u}(w) = E(C_{K_u}(w))$.
Since $C^\o_M(u)$ is perfect, $C^\o_M(u) \leq K_u$,
 so $C^\o_{K_u}(w) = C^\o_M(u)$.

By Table \ref{invSp4}, $E(C_M(u))/\gen{x} \cong \SL_2 * \SL_2$
 and the centralizer of an element of order four which acts by
 an involutive automorphism is $\SL_2 * k^*$.
As this does not match the centralizer of $w$,
 the action of $w$ is induced via conjugation by either $x$ or $v$.
So $E(C_M(u)) \cong L_1 \times L_2$ with $L_i \cong \SL_2$,
 $u\in Z(L_1)$, and $v\in L_2$, while $x$ belongs to neither $L_1$ nor $L_2$.
On the other hand, $E(C_M(u)) = E(C_{K_u}(w))$ can be viewed
 as either $E(C_{K_u}(x))$ or $E(C_{K_u}(v))$.
Since $u$ is the central involution of $K_u$,
 the same argument shows that $u$ belongs to neither $L_1$ nor $L_2$,
 a contradiction.
\end{analysis}

\begin{scase}\label{Lieranktwocomp_G2}
$\components^*$ contains a component $M \cong \G_2$.
\end{scase}

\begin{analysis}
By Table \ref{invG2},
 the three involutions in $\Omega_1(S^\o)$ are conjugate in $M$, and
$$ E(C_M(x)) \cong \SL_2 * \SL_2
  \quad\textrm{for every $x\in \Omega_1(S^\o)^\#$\mathperiod} $$
By $(\star)$,
 $E(C_G(z)) \not\leq M$ for some $z\in \Omega_1(S^\o)$.
Since $\pr(E(C_G(z))=2$,
 the group $E(C_G(z))$ must be quasisimple.
Since $S^\o \leq C^\o_G(z)$, we find $z\in Z(E(C_G(z)))$.
By Table \ref{findSp4}, $E(C_G(z)) \cong \Sp_4$.
Since $E(C_M(z))$ is perfect and
 $M\in \components$ induces descent for $A^\#$,
 $E(C_G(z)) \in \components^*$ too.
The result now follows by reduction to \caseref{Lieranktwocomp_Sp4}.
\end{analysis}

Before attacking the final two cases, we observe that pairs of components
of $\pPSL_2$ type always generate a proper subgroup.

\begin{lemma}\label{Lieranktwocomp_amalgam}
For any two components $L,J\in \components$ of $\pPSL_2$ type,
 one of the following holds.
\begin{conclusions}
\item $E(C_G(i)) = L * J$ for some $i\in A^\#$.
\item $L$ and $J$ are algebraic subgroups of $\gen{L,J} \in \components^*$.
\end{conclusions}
\end{lemma}

\begin{proof}
Since $A$ centralizes $T$,
 $A$ fixes all components of $E(C_G(y))$ for any $y\in A^\#$,
 and $A$ normalizes $L$ and $J$.
Since $L$ and $J$ have no graph automorphisms, 
 $A$ acts on $L$ and $J$ via inner automorphisms
 by Fact \ref{autalg}.
We observe that $T \cap L$ is a Sylow\o 2-subgroup of $L$
 which is centralized by $A$.
Since a four-group in $\PSL_2$ normalizing a 2-torus 
 has an involution inverting the torus and $A$ centralizes the 2-torus,
 we have that $m(C_A(L)) \geq m(A)-1 \geq 2$ and $m(C_A(J)) \geq 2$.
Thus $C_A(L) \cap C_A(J) \neq 1$ and
 $\gen{L,J} \leq C_G(i)$ for some involution $i \in C_A(L) \cap C_A(J)$.
Since $L$ and $J$ induce descent for $A^\#$,
 they must lie inside components of $C_G(i)$ coming from $\components$.

Now there is a $y\in A^\#$ such that $L \normal E(C_G(y))$.
Since $y$ normalizes the component $L' \in \components$
 with $L' \normal E(C_G(i))$ which contains $L$,
Fact \ref{autalg} says that $y$ acts on $L'$ by
 a graph automorphism composed with an inner automorphism.
So $C^\o_{L'}(y)$ is a reductive algebraic subgroup of $L'$
 by Fact \ref{Creductive2},
 which has $L$ as a normal subgroup, so $L$ is a component of $C^\o_{L'}(y)$.
If $[L,J]=1$ then $L = L'$ and $J = J'$ (defined similarly),
 so the first conclusion follows.
If $[L,J] \neq 1$ then $L'$ contains $J$ as an algebraic subgroup too,
 so the second conclusion follows.
\end{proof}

We now return to our case analysis.

\begin{scase}\label{Lieranktwocomp_PSp4}
$\components^*$ contains a component $M \cong \PSp_4$.
\end{scase}

\begin{analysis}
By Table \ref{invSp4}, the three involutions $x,y,z \in \Omega_1(S^\o)$
 have centralizers $C^\o_M(z) \cong \SL_2 * \SL_2$ and
 $C^\o_M(x) \cong C^\o_M(y) \cong \PSL_2 * k^*$
 with $x$ and $y$ conjugate in $M$.

Suppose first that $E(C_G(z)) \not\leq M$.
There is an $i\in A^\#$ with $M = E(C_G(i))$.
Since $M$ induces descent for $A^\#$, the two components
 of $C^\o_M(z) \cong \SL_2 * \SL_2$ are found in $\components$.
Since these two components induce descent for $A^\#$,
 $C^\o_M(z) \leq K_z := E(C_G(z))$ and
 the components of $K_z$ are found in $\components$.
Since $i$ normalizes $K_z$,
 the group $C_{K_z}(i)$ is an algebraic subgroup of $K_z$ by Fact \ref{autalg}.
So $C_{K_z}(i) \cong \SL_2 * \SL_2$ too, and $K_z$ is quasisimple.
By Table \ref{findSp4}, $E(C_G(z)) \cong \Sp_4$.
Now $E(C_G(z)) \in \components^*$,
 and the result follows by reduction to \caseref{Lieranktwocomp_Sp4}.
So $E(C_G(z)) \leq M$ and $E(C_G(z)) \cong \SL_2 * \SL_2$.

Therefore, $E(C_G(x)) \not\leq M$ by $(\star)$.
Now $E(C_G(x))$ must be either $\Sp_4$ or $\SL_2 \times \PSL_2$.
Since $E(C_M(x))$ is perfect and $M\in \components$ induces descent for $A^\#$,
 the $\Sp_4$ case reduces to \caseref{Lieranktwocomp_Sp4}.
So $E(C_G(x)) \cong \SL_2 \times \PSL_2$.

Let $J$ be the component of $E(C_G(x))$ which is isomorphic to $\SL_2$.
If $J$ commuted with either $\SL_2$ component of $E(C_G(z))$,
 then $J \leq E(C_G(z))$ since $[x,z]=1$.
Now $J$ does not commute with either $\SL_2$ component of
 $E(C_G(z))$, because otherwise it would be one of them.
Now let $L$ be a component of $E(C_G(z))$ with $L \cong \SL_2$.
By Lemma \ref{Lieranktwocomp_amalgam},
 $K := \gen{L,J} \in \components^*$.
In this group, the two involutions $x,z\in \Omega_1(S^\o)$ have $\SL_2$
 type components $J$ and $L$ in $C_G(x)$ and $C_G(z)$, respectively.
Thus $K \not\cong \PSp_4$ where only one toral involution may
 have a component of type $\SL_2$ by Table \ref{invSp4}.
We may also assume that $K$ is not isomorphic to either $\Sp_4$ or $\G_2$
 by reduction to Cases \ref{Lieranktwocomp_Sp4} and \ref{Lieranktwocomp_G2}
 since $L \in \components$ induces descent for $A^\#$.
So $K \cong \pPSL_3$.  Now the involutions of $\Omega_1(S^\o)$ are
$K$-conjugate by Table \ref{invSL3}, contradicting the fact that
$E(C_G(z)) \not\cong E(C_G(x))$.
\end{analysis}

Therefore we have only two possibilities.

\begin{scase}\label{Lieranktwocomp_pPSL3}
Every component $M \in \components^*$ is isomorphic to $\SL_3$ or $\PSL_3$.
\end{scase}

\begin{analysis}
Fix $M\in \components^*$.
Since $M \cong \pPSL_3$, the three distinct involutions $x,y,z\in
\Omega_1(S^\o)$ are conjugate in $M$.
For any one of these involutions $v\in \Omega_1(S^\o)$,
we define $L_v = E(C_M(v)) \cong \SL_2$.
Since $M \cong \pPSL_3$, we have $v\in Z(E(C_M(v)))$ and
 $L_u \neq L_v$ for $u \neq v$ (see Table \ref{invSL3}).
Since $E(C_G(v)) \not\cong \Sp_4$ by \caseref{Lieranktwocomp_Sp4},
 we have $L_v \normal E(C_G(v))$ and
 there is a second component in $E(C_G(v))$, which we denote by $J_v$.
Since $L_v \leq E(C_M(\Omega_1(S^\o) \cap J_v))$,
 we find $\Omega_1(S^\o) \cap J_v = \gen{v}$, and
 thus $J_u \neq J_v$ for $u \neq v$.
Since $J_x,J_y,J_z \not\leq M$,
 the set $\mathcal L = \{L_x,J_x,L_y,J_y,L_z,J_y\}$
 of these subgroups has exactly 6 elements.
We also see that $J_v \cong \SL_2$ and
 $E(C_G(v)) \cong \SL_2 * \SL_2$ for $v\in \Omega_1(S^\o)^\#$.

Let $H := C^\o_G(S^\o)$.
Since $H$ centralizes $\Omega_1(S^\o)$,
 $H$ normalizes any $L\in \mathcal L$,
 $H$ normalizes $M = \gen{L_x,L_y,L_z}$.
Since $H$ is connected and definable,
 $H = C^\o_H(M) (M \cap H)^\o$ by Fact \ref{autalg}.
So for any component $K\in \components^*$,
 there is a natural embedding of the Weyl group
 $W(K) := N_K(C_K(S^\o)) / C_K(S^\o)$ of $K$ into $N_G(H) / H$.
We define the ``Weyl group'' $W$ of $G$ to be the subgroup of $N_G(H)/H$
 which is generated by the Weyl groups $N_K(H)/H$ of each component
 $K\in \components^*$.
We observe that the subgroups $L \in \mathcal L$ are normalized by $H$,
 and therefore are the root $\SL_2$-subgroups of those components
 $K\in \components^*$ to which they belong.
Since $K \cong \pPSL_3$ by previous cases,
 $W(K)$ is generated by $N_{H L}(H)/H$ for some $L\in \mathcal L$,
 so $W$ is generated by $N_{H L}(H)/H$ for $L\in \mathcal L$.
Now $|N_{H L}(H)/H| = |N_L(H \cap L)/(H \cap L)| = 2$. 
For any $L\in \mathcal L$, we let $r_L$ denote the involution of
$N_{H L}(H)/H$.  We observe that $r_L$ acts on $S^\o$ by reflection
 in the sense that $S^\o = C_{S^\o}(r_L) \times [S^\o,r_L]$
 with $\pr(C_{S^\o}(r_L)) = 1$ and $\pr([S^\o,r_L]) = 1$.

By the Tate module argument \cite[\S3.3]{Be01,BB01}
 (see also \cite[\S2.4]{BBBC07}),
 there is a (not necessarily faithful) representation of $W$
 by $2\times2$ matrices over the 2-adic integers $\Z_2$,
 and over $\C$ after tensoring.
We let $\ov W$ denote the image of $W$ inside $\End(\C^2)$.
The involutions $r_L$ act on $\C^2$ as reflections
 with one $+1$ eigenvalue and one $-1$ eigenvalue.
The kernel $U$ of the natural action of $W$ on $\mathcal L$ leaves invariant
 the three subgroups $S^\o \cap L_x$, $S^\o \cap L_y$, and $S^\o \cap L_z$,
 and hence elements of $U$ act either trivially or by inverting $S^\o$.
So the action of $W$ on $\mathcal L$ can be factored through $\ov W$,
 and $\ov U$ consists of the scalar matrices $\pm1$.

We now consider the Weyl group $W(M) = N_M(C_M(S^\o)) / C_M(S^\o)$ of $M$.
Let $r\in N_M(C_M(S^\o)) \cap L_x$ be a representative for the Weyl group
element of $L_x$ associated to the maximal torus $C_M(S^\o)$.
We may assume $r\in S$ by \cite[Ex.~11 p.~93]{BN}.
Since $C_M(x) \cong \SL_2 * k^*$ with amalgamation by Table \ref{invSL3},
 conjugation by $r$ must swap $y$ and $z$.
So $[L_x,J_y] \neq 1$.
By Lemma \ref{Lieranktwocomp_amalgam},
 $\gen{L_x,J_y} \in \components^*$.
Similarly, $\gen{L_y,J_z},\gen{L_z,J_x} \in \components^*$.
For any $K\in \components^*$,
 we have $K \cong \pPSL_3$ by the previous cases.
Since $\mathcal L$ is independent of the choice of $M$,
 $K$ has the form $\gen{L,L'}$ for two $L,L'\in \mathcal L$.
There is no $\SL_2 * \SL_2$ contained in the centralizer of an involution
 in $K \cong \pPSL_3$ by Table \ref{invSL3}, so
 $\gen{L_v,J_v} \notin \components^*$ for any $v\in \Omega_1(S^\o)^\#$.
For any $u,v,w\in \Omega_1(S^\o)^\#$ with $u \neq v$,
the group $J_w$ is not contained in $M = \gen{L_u,L_v}$.

Assembling these various facts, we discover that 
$$ \components^* =
  \{ \gen{L_x,L_y,L_z}, \gen{L_x,J_y,J_z}, 
     \gen{J_x,L_y,J_z}, \gen{J_x,J_y,L_z} \}\mathperiod $$
Thus we have the following geometry.
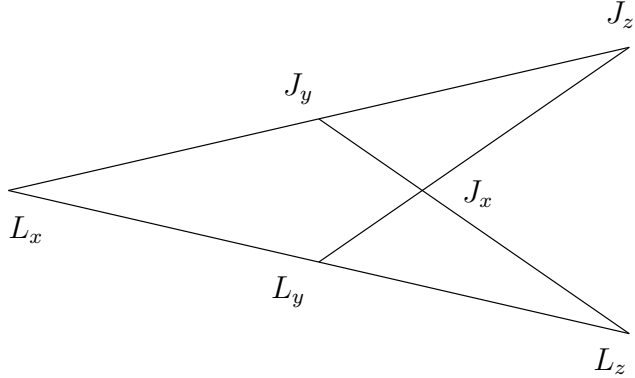
\begin{center}
\providecommand\color[2][]{}
\begin{picture}(0,0)%
\epsfig{file=Sym4.pstex}%
\end{picture}%
\setlength{\unitlength}{3947sp}%
\begingroup\makeatletter\ifx\SetFigFont\undefined%
\gdef\SetFigFont#1#2#3#4#5{%
  \reset@font\fontsize{#1}{#2pt}%
  \fontfamily{#3}\fontseries{#4}\fontshape{#5}%
  \selectfont}%
\fi\endgroup%
\begin{picture}(4208,2389)(5989,-3044)
\put(8851,-1936){\makebox(0,0)[lb]{\smash{{\SetFigFont{12}{14.4}{\rmdefault}{\mddefault}{\updefault}{\color[rgb]{0,0,0}$J_x$}%
}}}}
\put(6001,-2161){\makebox(0,0)[lb]{\smash{{\SetFigFont{12}{14.4}{\rmdefault}{\mddefault}{\updefault}{\color[rgb]{0,0,0}$L_x$}%
}}}}
\put(9676,-2986){\makebox(0,0)[lb]{\smash{{\SetFigFont{12}{14.4}{\rmdefault}{\mddefault}{\updefault}{\color[rgb]{0,0,0}$L_z$}%
}}}}
\put(9751,-811){\makebox(0,0)[lb]{\smash{{\SetFigFont{12}{14.4}{\rmdefault}{\mddefault}{\updefault}{\color[rgb]{0,0,0}$J_z$}%
}}}}
\put(7726,-1261){\makebox(0,0)[lb]{\smash{{\SetFigFont{12}{14.4}{\rmdefault}{\mddefault}{\updefault}{\color[rgb]{0,0,0}$J_y$}%
}}}}
\put(7651,-2536){\makebox(0,0)[lb]{\smash{{\SetFigFont{12}{14.4}{\rmdefault}{\mddefault}{\updefault}{\color[rgb]{0,0,0}$L_y$}%
}}}}
\end{picture}%

\end{center}

For each $K\in \components^*$, the group $N_K(H)/H$ acts 2-transitively
on the three subgroups from $\mathcal L$ which are contained in $K$.
Recall that $U$ is the kernel of the action of $W$ on $\mathcal L$.
So the permutation group $W/U$ preserves the geometry consisting of
the six point in $\mathcal L$ and the four lines in $\components^*$,
and permutes 2-transitively the points on each line.
Since a permutation of the lines $\components^*$ determines a permutation
of the points $\mathcal L$, we have an injective homomorphism
$$ \bar{W}/\bar{U} \hookrightarrow \Sym_{\components^*} \cong \Sym_4\mathperiod $$
Since $W/U$ acts 2-transitively on the points on each line,
 $\bar{W}/\bar{U} \cong \Sym_4$.

Let $\bar{W}_0 \leq \bar{W}$ be the preimage of $\Alt_4$ under the
quotient map.  We recall that $\bar{U}$ consists of scalar $\pm1$ matrices.
If $\bar{U} = 1$ then $\bar{W}_0 \cong \Alt_4$ is a subgroup of
$\GL_2(\C)$, so the subgroup $\Z/4\Z \leq \Alt_4$ is diagonalizable and
contains the scalar $-1$ matrix, a contradiction because $\Alt_4$ is centerless.
Thus $\bar{U} \neq 1$.

Let $\bar{Q}$ be the pull-back of the group $(\Z/2\Z)^2 \normal \Alt_4$
 under the quotient map $\bar{W}_0 \to \bar{W_0/\bar{U}}$.
Then $\bar{Q}$ is acted upon by $\Z/3\Z$, and
 the elements of $\bar{Q} \setminus \bar{U}$ have order 2 or 4.
If the elements have order 2, then $\bar{Q} \cong (\Z/2\Z)^3$,
which is not a subgroup of $\GL_2(\C)$.
So these elements have order 4, and $\bar{Q} \cong Q_8$. 
Now $\bar{W}_0$ is the semidirect product $Q_8 \rtimes (\Z/3\Z)$ of
 the quaternion group $Q_8$ of order 8 and the cyclic group of order 3.
We also find that the subgroup $\bar{Q} \cong Q_8$ from $\bar{W}_0$
 is normal in $\bar{W}$ and $\bar{W} / \bar{Q} \cong \Sym_3$.

Now the group $\bar{W}$ has order $2 \cdot 4! = 48$, contains
 a conjugacy class of reflections $r_L$ for $L\in \mathcal L$,
 and has $\abs{Z(\bar{W})} = 2$.
Such a group does not exist, as we now show in a couple ways.

By the classification of irreducible complex reflection groups
\cite[Table VII, p.~301]{ShTo54}, 
 $\bar{W}$ must be one of $G(12,6,2)$, $G(24,24,2)$,
 or have Shephard--Todd number either 6 or 12,
 as only these have order 48.
We can determine from \cite[Table VII\ p.\ 301]{ShTo54} that
$G(12,6,2)$ and number 6 have centers or order 2, while $G(24,24,2)$
 and number 12 have no conjugacy class of exactly six reflections.
One can preform these computations with the following
GAP \cite{GAP4} commands.
\begin{verbatim}
RequirePackage("chevie");
Size(Centre(ComplexReflectionGroup(12,6,2)));
Size(Centre(ComplexReflectionGroup(6)));
W := ComplexReflectionGroup(24,24,2);
ForAny(Reflections(W), x -> Size(ConjugacyClass(W,x))=6);
W := ComplexReflectionGroup(12);
ForAny(Reflections(W), x -> Size(ConjugacyClass(W,x))=6);
\end{verbatim}

To check the result by hand as follows.
There are six such reflections, so $\abs{C_{\bar{W}}(r_L)} = 8$,
$\abs{C_{\bar{Q}}(r_L)} = 4$, and $C_{\bar{Q}}(r_L) \cong \Z/4\Z$.
We can choose a basis in the Tate module $V$ so that $r_L$ is
 represented by $\diag(-1,1)$.
So $C_{\bar{Q}}(r_L)$ consists of diagonal matrices too, and
 must be generated by either $t = \diag(-i,i)$ or $\diag(i,-i)$.
So $r_L t$ is a scalar matrix $\diag(\pm1,\pm1)$ and belongs to $Z(\bar{W})$.
This means $r_L \in Z(\bar{W}) \bar{Q} \normal \bar{W}$ and
 the conjugates of $r_L$ can not generate the group, a contradiction.
\end{analysis}

This concludes the proof of the Trichotomy Theorem.

\small
\bibliographystyle{amsplain}
\bibliography{burdges,fMr}

\def\cprime{$'$}
\providecommand{\bysame}{\leavevmode\hbox to3em{\hrulefill}\thinspace}
\providecommand{\MR}{\relax\ifhmode\unskip\space\fi MR }
\providecommand{\MRhref}[2]{%
  \href{http://www.ams.org/mathscinet-getitem?mr=#1}{#2}
}
\providecommand{\href}[2]{#2}
\begin{thebibliography}{10}

\bibitem{ABC97}
Tuna Alt{\i}nel, Alexandre Borovik, and Gregory Cherlin, \emph{Groups of mixed
  type}, J. Algebra \textbf{192} (1997), no.~2, 524--571. \MR{98d:03047}

\bibitem{Asch}
Michael Aschbacher, \emph{Finite group theory}, Cambridge Studies in Advanced
  Mathematics, vol.~10, Cambridge University Press, Cambridge, 1993, Corrected
  reprint of the 1986 original. \MR{95b:20002}

\bibitem{Bel87}
Oleg~V. Belegradek, \emph{On groups of finite {M}orley rank}, Abstracts of the
  Eight International Congress of Logic, Methodology and Philosophy of Science
  LMPS'87 (Moscow), 1987, 17--22 August 1987, pp.~100--102.

\bibitem{Be01}
Ay{\c{s}}e Berkman, \emph{The classical involution theorem for groups of finite
  {M}orley rank}, J. Algebra \textbf{243} (2001), no.~2, 361--384.
  \MR{2002e:20060}

\bibitem{BBBC07}
Ay{\c{s}}e Berkman, Alexandre Borovik, Jeffrey Burdges, and Gregory Cherlin,
  \emph{A generic identification theorem for ${L}^*$-groups of finite {M}orley
  rank}, J. Algebra (2007), \Toappear.

\bibitem{BB01}
Ay{\c{s}}e Berkman and Alexandre~V. Borovik, \emph{A generic identification
  theorem for groups of finite {M}orley rank}, J. London Math. Soc. (2)
  \textbf{69} (2004), no.~1, 14--26. \MR{MR2025324}

\bibitem{BP}
Aleksandr~Vasilievich Borovik and Bruno~Petrovich Poizat, \emph{Tores et
  {$p$}-groupes}, J. Symbolic Logic \textbf{55} (1990), no.~2, 478--491.
  \MR{91j:03045}

\bibitem{Bo95}
Alexandre Borovik, \emph{Simple locally finite groups of finite {M}orley rank
  and odd type}, Finite and locally finite groups (Istanbul, 1994), NATO Adv.
  Sci. Inst. Ser. C Math. Phys. Sci., vol. 471, Kluwer Acad. Publ., Dordrecht,
  1995, pp.~247--284. \MR{96h:20061}

\bibitem{BBC}
Alexandre Borovik, Jeffrey Burdges, and Gregory Cherlin, \emph{Involutions in
  groups of finite {M}orley rank of degenerate type}, Selecta (2007),
  \Toappear.

\bibitem{BBN04}
Alexandre Borovik, Jeffrey Burdges, and Ali Nesin, \emph{Uniqueness cases in
  odd type groups of finite {M}orley rank}, J. London Math. Soc. (2006),
  \Toappear.

\bibitem{BN}
Alexandre Borovik and Ali Nesin, \emph{Groups of finite {M}orley rank}, The
  Clarendon Press Oxford University Press, New York, 1994, Oxford Science
  Publications. \MR{96c:20004}

\bibitem{BuPhd}
Jeff Burdges, \emph{Simple groups of finite {M}orley rank of odd and degenerate
  type}, Ph.D. thesis, Rutgers University, New Brunswick, New Jersey, 2004.

\bibitem{Bu03}
Jeffrey Burdges, \emph{A signalizer functor theorem for groups of finite
  {M}orley rank}, J. Algebra \textbf{274} (2004), no.~1, 215--229.

\bibitem{Bu07a}
\bysame, \emph{Signalizers and balance in groups of finite {M}orley rank},
  \Submitted[J. LMS], 2007.

\bibitem{BCJ}
Jeffrey Burdges, Gregory Cherlin, and Eric Jaligot, \emph{Minimal connected
  simple groups of finite {M}orley rank with strongly embedded subgroups}, J.
  Algebra (2007), \Toappear.

\bibitem{Carter93}
Roger~W. Carter, \emph{Finite groups of {L}ie type}, Wiley Classics Library,
  John Wiley \& Sons Ltd., Chichester, 1993, Conjugacy classes and complex
  characters, Reprint of the 1985 original, A Wiley-Interscience Publication.
  \MR{94k:20020}

\bibitem{CJ01}
Gregory Cherlin and Eric Jaligot, \emph{Tame minimal simple groups of finite
  {M}orley rank}, J. Algebra \textbf{276} (2004), no.~1, 13--79.

\bibitem{Cohen76}
Arjeh~M. Cohen, \emph{Finite complex reflection groups}, Ann. Sci. \'Ecole
  Norm. Sup. (4) \textbf{9} (1976), no.~3, 379--436. \MR{54 \#10437}

\bibitem{GAP4}
The GAP~Group, \emph{{GAP -- Groups, Algorithms, and Programming, Version
  4.3}}, 2002, \verb+(http://www.gap-system.org)+.

\bibitem{GLS3}
Daniel Gorenstein, Richard Lyons, and Ronald Solomon, \emph{The classification
  of the finite simple groups. {N}umber 3. {P}art {I}. {C}hapter {A}},
  Mathematical Surveys and Monographs, vol.~40, American Mathematical Society,
  Providence, RI, 1998, Almost simple $K$-groups. \MR{98j:20011}

\bibitem{Ne91}
Ali Nesin, \emph{Generalized {F}itting subgroup of a group of finite {M}orley
  rank}, J. Symbolic Logic \textbf{56} (1991), no.~4, 1391--1399.
  \MR{92h:03049}

\bibitem{ShTo54}
G.~C. Shephard and J.~A. Todd, \emph{Finite unitary reflection groups},
  Canadian J. Math. \textbf{6} (1954), 274--304. \MR{15,600b}

\bibitem{St2}
Robert Steinberg, \emph{Endomorphisms of linear algebraic groups}, Memoirs of
  the American Mathematical Society, No. 80, American Mathematical Society,
  Providence, R.I., 1968. \MR{37 \#6288}

\bibitem{WalterClsInv}
John~H. Walter, \emph{The {$B$}-conjecture; characterization of {C}hevalley
  groups}, Mem. Amer. Math. Soc. \textbf{61} (1986), no.~345, iv+196.
  \MR{87h:20035}

\end{thebibliography}

\affiliationone{
Alexandre Borovik and Jeffrey Burdges\\
School of Mathematics, The University of Manchester\\
PO Box 88, Sackville St., Manchester M60 1QD, England\\
\email{Alexandre.Borovik@manchester.ac.uk\\ Jeffrey.Burdges@manchester.ac.uk}}

\end{document}